\documentstyle[12pt]{article}
\begin{document}
\section*{ Introduction}
\thispagestyle{empty}
  Vector bundles on complex surfaces have been
extensively studied by means of several different methods.
See for example the books of
Kobayashi [12] and Okonek, Schneider, Spindler [13].
Stable holomorphic bundles on a K\"ahler surface
correspond by a theorem of Donaldson [5]
to irreducible anti-self-dual connections on the surface.
This result connects the study of holomorphic
vector bundles with moduli space of instantons.
{}From the point of view of the study of instantons,
vector bundles on the blow up of ${\bf P}^2$
appear in Hurtubise's paper on Instantons and Jumping Lines
[10] and in Boyer-Hurtubise-Milgram-Mann [1] in
their proof of the Atiyah-Jones conjecture.
Subsequently Hurtubise and Milgram [11]
proved an extended version of the
 Atiyah-Jones conjecture for ruled surfaces, by means of studying the
structure of holomorphic bundles on ruled
surfaces.
\\
To motivate our study of holomorphic vector bundles
on blow-ups from a different stand point
we mention a fundamental result on the classification
of rational surfaces, see Griffiths and Harris [8].
\\
{ \bf Theorem} : Every rational surface is obtained by blowing up
points on either ${\bf P}^2$ or on a rational ruled
surface.
\\
The previous theorem
suggests that
the understanding
of vector bundles on rational surfaces
depends  on the analysis of the behavior of
vector bundles under blow-ups.
A large amount of work has been done on  vector
bundles on ${\bf P}^2$ ( see for example
the book by  Okonek, Schneider, Spindler [13]).
In a sense we  can also say that vector bundles on ruled surfaces
are well understood ( see Brosius [2] [3],
Qin [14], Hurtubise and  Milgram [11]).
Some examples of work on moduli spaces of
holomorphic vector bundles on blow-ups
are the papers by  by Freedman and Morgan [6][7],
Brussee [4], and Qin [15].
\\
The blow-up of a point on a surface is a local
operation  in the sense that one blows-up the point
inside one of its coordinate neighborhoods.
Roughly speaking we may see  the ``difference''
between moduli spaces of  bundles
on a rational surface and moduli spaces of
bundles on  one of its minimal models by studying bundles
on the blow up of ${\bf C}^2.$
\\
In this work we concentrate on the study
of bundles on blow-ups in the local sense,
that is in a neighborhood of the exceptional divisor.
Our approach is  quite concrete, as we
give bundles explictly by their transition matrices
and present the moduli spaces as quotients of
a vector space ${\bf C}^n$ by an equivalence relation.
\\
In Section 3 we construct a canonical form of transition matrix
for rank two bundles on the blow up of  ${\bf C}^2.$ Namely, we prove
the following:
\\
\vspace {5 mm}
\noindent{\bf Theorem 2.1}:  Let $E$ be a
 holomorphic rank two vector  bundle on $ \widetilde{\bf C}^2  $ with
zero first Chern class and
let j be the integer that satisfies
$E_{\ell} \simeq {\cal O}(j) \oplus  {\cal O}(-j). $
(Where $E_{\ell}$ is the restriction of $E$ to
the exceptional divisor.)
Then  $E$  has a transition matrix
of the form
\\
$$\left(\matrix {z^j & p \cr 0 &  z^{-j} \cr }\right)$$
from $U$ to $V,$  where
\\
$$p = \sum_{i = 1}^{2j-2} \sum_{l = i-j+1}^{j-1}p_{il}z^lu^i.$$
In particular $p$ depends on a finite number of parameters.
\\
\vspace{5 mm}
\\
We then define the moduli space ${\cal M}_j$ as
the space of equivalence classes of such bundles having
restriction  ${\cal O}(j) \oplus  {\cal O}(-j)$ to the
exceptional divisor, modulo holomorphic equivalence.
It follows immediately from our canonical form
of a  transition matrix that:
\\
\vspace{5 mm}
\noindent{\bf Corollary 2.3}:  ${\cal M}_0$  consists of a single point.
\\
\vspace{5 mm}
\\
\noindent{\bf Corollary 2.5}:  ${\cal M}_1$  consists of a single point.
\\
\vspace{5 mm}
\\
In Section 4 we continue the study of  ${\cal M}_j$
for $j \ge 2.$
To do this we analyze the problem of when two holomorphic
bundles in  ${\cal M}_j$ are isomorphic.a simple  characterization on
the first formal neighborhood
of the exceptional divisor.
\\
\vspace{ 5 mm}
\noindent{\bf Proposition 3.3}
On the first formal neighborhood,
two holomorphic bundles $E^{(1)}$ and $E^{(1)\prime}$ with
transition matrices
$$\left(\matrix{z^j & p_1   \cr 0 & z^{-j} \cr}\right)$$
and
$$\left(\matrix{z^j & p^\prime_1   \cr 0 & z^{-j} \cr}\right)$$
 respectively  are isomorphic iff
$p^\prime_1 = \lambda p_1$ for some
 $\lambda \in  {\bf C} - \{0\}$.
\\
\vspace {5 mm}
\\
Once one passes the first formal neighborhood,
the holomorphic equivalences become more intricate.
In 5.1 we give a detailed description of
${\cal M}_2.$ Topologicaly, we have:
\\
\vspace{5 mm}
\\
\noindent{\bf Theorem 4.2}: The moduli space ${\cal M}_2$ is
homeomorphic to
 the union ${\bf P}^1 \cup \{p,q\},$
of a complex projective plane ${\bf P}^1$ and two points
  with  a basis  of open sets given by
$${\cal U}  \cup \{p,U : U \in {\cal U} - \phi \} \cup
\{p,q,U  : U \in {\cal U} - \phi \}  $$
where ${\cal U}$ is a basis for the standard topology on  ${\bf P}^1.$
\\
\vspace{5 mm}
\\
In 5.2 we describe ${\cal M}_3$ and in 5.3 we give the generic
description of ${\cal M}_j.$ Our general results are:
\\
\vspace {5 mm}
\noindent{\bf Theorem 4.4}
The generic set of the moduli space ${\cal M}_j$ is a complex
projective space of dimension $2j-3$ minus a closed
subvariety of complex codimension bigger than or equal to two.
\\
\vspace{5 mm}
\\
\noindent{\bf Remark 4.6}: The moduli space ${\cal M}_j$ also contains
complex projective spaces of every dimension smaller
than $2j-3, $ each minus some closed subvariety.
\\
\vspace {5 mm}
\\
\noindent{\bf Remark 4.7}: If we give  ${\cal M}_j$ the topology
induced
from ${\bf C}^N$, then  ${\cal M}_j$ is not a Hausdorff space.
For example, the direct sum bundle given by
$\left(\matrix{ z^j & 0 \cr 0 & z^{-j} \cr
}\right)$ is arbitrarily close to any other bundle.
\\
\vspace {5 mm}
\\
\noindent{\bf Remark 4.8}: Note that the word generic here is used
in the sense that the moduli space  ${\cal M}_j$
consists of subsets out of
which ${\bf P}^{2j-3} $ is the subset of highest
dimension.
\\
\vspace{5mm}
\\
Finally  in Section 6 we give  some examples
of the result of building up bundles on the
blow up of a compact surface using our
canonical form of a transition matrix for a neighborhood of the
exceptional divisor.
\\
\vspace{5mm}
\\
Note: This is a quite long file, so I am only sending the "introduction."
If anyone wants the whole file, be welcome to write to
gasparim@ictp.trieste.it
\\
\end{document}